\documentclass[12pt]{article}
\pdfoutput=1
\usepackage{amsmath}    
\usepackage{graphicx}   
\usepackage{verbatim}   
\usepackage{color}      
\usepackage{subfigure}  
\usepackage{amssymb}
\usepackage{cite}
\usepackage{scrextend}
\usepackage[english]{babel}
\usepackage{blindtext}
\newtheorem{theorem}{Theorem}[section]
\newtheorem{proposition}[theorem]{Proposition}

\newtheorem{corollary}[theorem]{Corollary}
\newtheorem{lemma}[theorem]{Lemma}

\newcommand{\proof}{\noindent{\bf Proof.\ }}
\newcommand{\qed}{\hfill $\square$\medskip}

\DeclareMathOperator {\gp} {gp}
\DeclareMathOperator {\ip} {ip}
\DeclareMathOperator {\ic} {ic}
\DeclareMathOperator {\diam} {diam}

\begin{document}
	
\title{Graph theory general position problem}
	
\author{
	Paul Manuel $^{a}$
	\and
	Sandi Klav\v zar $^{b,c,d}$
}

\date{}

\maketitle
\vspace{-0.8 cm}
\begin{center}
	$^a$ Department of Information Science, College of Computing Science and Engineering, Kuwait University, Kuwait \\
	{\tt pauldmanuel@gmail.com}\\
	\medskip
	
	$^b$ Faculty of Mathematics and Physics, University of Ljubljana, Slovenia\\
	{\tt sandi.klavzar@fmf.uni-lj.si}\\
	\medskip
	
	$^c$ Faculty of Natural Sciences and Mathematics, University of Maribor, Slovenia\\
	\medskip
	
	$^d$ Institute of Mathematics, Physics and Mechanics, Ljubljana, Slovenia\\
	
\end{center}

\begin{abstract}
The classical no-three-in-line problem is to find the maximum number of points that can be placed in the $n \times n$ grid so that no three points lie on a line. Given a set $S$ of points in an Euclidean plane, the General Position Subset Selection Problem is to find a maximum subset $S'$ of $S$ such that no three points of $S'$ are collinear. Motivated by these problems, the following graph theory variation is introduced: Given a graph $G$, determine a largest set $S$ of vertices of $G$ such that no three vertices of $S$ lie on a common geodesic. Such a set is a gp-set of $G$ and its size is the gp-number $\gp(G)$ of $G$. Upper bounds on $\gp(G)$ in terms of different isometric covers are given and used to determine the gp-number of several classes of graphs. Connections between general position sets and packings are investigated and used to give lower bounds on the gp-number. It is also proved that the general position problem is NP-complete. 
\end{abstract}


\medskip
\noindent{\bf AMS Subj.\ Class.: 05C12, 05C70, 68Q25} 

\section{Introduction}
\label{sec:introduction}

The no-three-in-line problem is to find the maximum number of points that can be placed in the $n \times n$ grid so that no three points lie on a line. This celebrated century-old problem that was posed by Dudeney~\cite{dudeney-1917} is still open. For some recent related developments, see~\cite{PoWo07, misiak-2016} and references therein. In the first of these two papers the problem is extended to 3D, while in the second it is proved that at most $2{\rm gcd}(m,n)$ points could be placed with no three in a line on an $m \times n$ discrete torus. The no-three-in-line problem was in discrete geometry extended to the General Position Subset Selection Problem~\cite{FrKa16, PaWo13}, where for a given set of points in the plane one aims to determine a largest subset of points in general position. In~\cite{FrKa16} it is proved, among other results, that the problem is NP- and APX-hard, while in~\cite{PaWo13} asymptotic bounds on the function $f(n,\ell)$ are derived, where $f(n,\ell)$ is the maximum integer such that every set of $n$ points in the plane with no more than $\ell$ collinear contains a subset of $f(n,\ell)$ points with no three collinear. 

The above problems motivated us to define a similar problem in graph theory as follows: Given a graph $G$, the {\em graph theory general position problem} is to find a largest set of vertices $S\subseteq V(G)$, such that no three vertices of $S$ lie on a common geodesic in $G$. Note that an intrinsic difference between the discrete geometry problem and the graph theory general position problem is that in the first case for given points $x$ and $y$ there is only one straight line passing through $x$ and $y$, while in the graph theory problem there can be several geodesics passing through two vertices. 

\begin{figure}[ht!]
	\begin{center}
		\scalebox{0.35}{\includegraphics{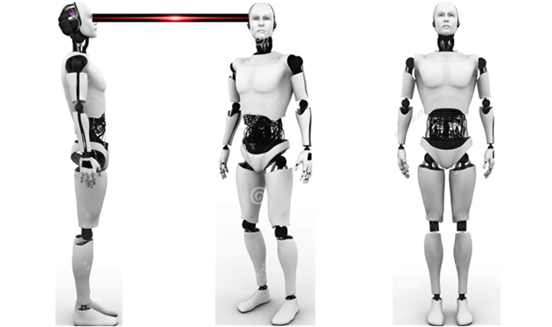}}
	\end{center}
	\caption{When three robots are in a straight line, the robot which stands in the middle blocks the sensor ray between the other two robots.}
	\label{fig:Robot_Application}
\end{figure}

Here is another motivation for the graph theory general position problem. Autonomous robots are intelligent machines that use sensors for processing visual signals and navigating in their environment~\cite{MaBaIn14}. Robots effectively see its neighbors using sensor systems~\cite[Chapter 5]{Corr16} such as infrared or ultrasound sensors. Navigation of a robot can be studied in a graph-structured framework~\cite{BaYeFe08, BaSe13}. The navigating agent can be assumed to be a point robot which moves from node to node of a ``graph space". Since the laser ray travels in a straight line, a robot $A$ cannot detect a robot $C$ when a third robot $B$ stands between $A$ and $C$ in the same line, see Fig.~\ref{fig:Robot_Application}. The controller needs to know how many robots can see each other. Given a collection of robots in a predefined structure, the problem is to find the largest number of robots who see each other which is just the graph theory general position problem. 

We proceed as follows. In the next section we give necessary definitions, general properties of general position sets, and exact values for the gp-number of some classes of graphs. In Section~\ref{sec:upper-bounds} upper bounds on the gp-number in terms of different isometric covers are obtained. It is also proved that the set of simplicial vertices of a block graph forms a maximum general position set. Relating general position sets with the diameter and the $k$-packing number, we derive in the subsequent section lower bounds on the gp-number. Then, in Section~\ref{sec:NPComp}, we prove that the general position problem is NP-complete.

\section{Preliminaries and examples}
\label{sec:preliminaries}

In this section we first define concepts and introduce the notation needed. Then we proceed to give the general position number of some families of graphs and along the way give some related general properties.  

All graphs considered in this paper are connected. The {\em distance} $d_G(u,v)$ between vertices $u$ and $v$ of a graph $G$ is the number of edges on a shortest $u,v$-path. Shortest paths are also known as {\em geodesics} or {\em isometric paths}. The {\em diameter} $\diam(G)$ of $G$ is the maximum distance between all pairs of vertices of $G$. A subgraph $H=(V(H),E(H))$ of a graph $G=(V(G),E(G))$ is {\em isometric} if $d_H(x,y) = d_G(x,y)$ holds for every pair of vertices $x,y$ of $H$. This is one of the key concepts in metric graph theory, cf.~\cite{beaudou-2009, pegrimek-2016, polat-2002, shpectorov-1998}. A {\em block} of a graph $G$ is a maximal connected subgraph of G that has no cut-vertex. A graph is a {\em block graph} if every block of it is complete. A vertex of a graph is {\em simplicial} if its neighbors induce a complete subgraph. For $n\in \mathbb{N}$ we will use the notation $[n] = \{1,\ldots,n\}$. 

A set $S$ of vertices of a graph $G$ is a {\em general position set} if no three vertices of $S$ lie on a common geodesic. A general position set $S$ of maximum cardinality is called a {\em $\gp$-set} of $G$. The cardinality of a gp-set of $G$ is called the {\em general position number} ({\em $\gp$-number} for short) of $G$ and denoted by $\gp(G)$.  

\begin{figure}[ht!]
	\begin{center}
		\scalebox{0.5}{\includegraphics{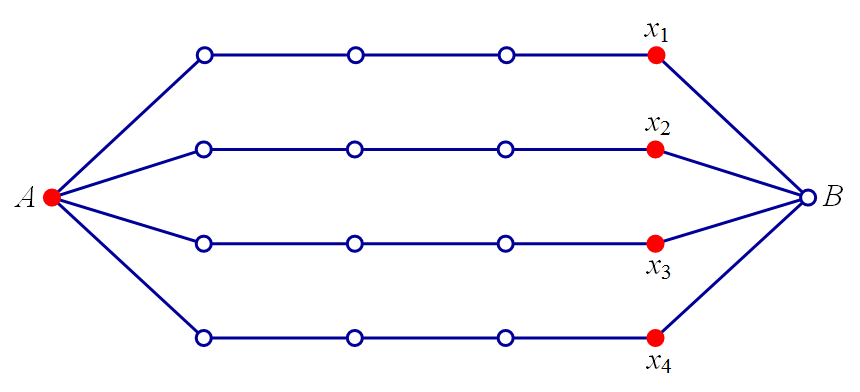}}
	\end{center}
	\caption{The theta graph $\Theta(4,5)$ and its gp-set.}
	\label{fig:theta-graph}
\end{figure}

As soon as $G$ has two vertices, $\gp(G)\ge 2$. For complete graphs, $\gp(K_n) = n$ for $n\ge 1$. Note also that $\gp(P_n) = 2$ for $n\ge 2$.  Consider next the cycle $C_n$ on vertices $v_1,\ldots, v_n$ with natural adjacencies. Let $S$ be an arbitrary general position set of $C_n$ and assume without loss of generality that $v_1 \in S$. Then 
\begin{quote}
$|S \cap \{v_2,v_3, \ldots, v_{\lceil(n+1)/2\rceil}\}| \leq 1$ and \\
$|S \cap \{v_{\lceil(n+1)/2\rceil+1}, v_{\lceil(n+1)/2\rceil+2},\ldots, v_n\}| \leq 1$. 
\end{quote}
If follows that $\gp(C_n)\leq 3$. If $n\ge 5$, then it is easy to find a gp-set in $C_n$ of order $3$. Hence $\gp(C_n) = 3$ for $n\ge 5$. Note also that $\gp(C_3) = 3$ and $\gp(C_4) = 2$. For $k\ge 2$ and $\ell\ge 2$, let $\Theta(k,\ell)$ be the graph consisting of two vertices $A$ and $B$ which are joined by $k$ internally disjoint paths each of length $\ell$. The vertices other than $A$ and $B$ are called internal vertices of $\Theta(k,\ell)$. See Fig.~\ref{fig:theta-graph} where $\Theta(4,5)$ is drawn. These graphs are known as {\em theta graphs}. 

\begin{proposition}
\label{prp:theta}
If $k\ge 2$ and $\ell\ge 3$, then $\gp(\Theta(k,\ell)) = k+1$. 
\end{proposition}

\proof
Let $R$ be a general position set of $\Theta(k, \ell)$. Let $P_i$, $i\in [k]$, denote the paths of $\Theta(k, \ell)$ joining $A$ and $B$. Consider arbitrary paths $P_i$ and $P_j$ of $\Theta(k, \ell)$. Then the union of $P_i$ and $P_j$ induces an isometric cycle $C$ of $\Theta(k, \ell)$ and hence $|R \cap V(C)| \le 3$. Therefore, if $R$ contains either $A$ or $B$, then each $P_i$, $i\in [k]$, contains at most one vertex of $R$ other than $A$ or $B$ respectively. If $R$ contains neither $A$ nor $B$, then only one path $P_j$ can contain two vertices from $R$ (clearly, it cannot contain three or more), and all the other paths $P_i$, where $i \neq j$, can have at most one vertex of $R$. In either case, $|R| \leq k+1$. Since $R$ is an arbitrary general position set of $\Theta(k, \ell)$, $\gp(\Theta(k,\ell)) \leq k+1$.
	
Let $x_1,\ldots, x_k$ be the vertices of $\Theta(k,\ell)$ that are adjacent to $B$ and set $S = \{A, x_1, \ldots, x_k\}$. See Fig.~\ref{fig:theta-graph}. It is easy to verify that $S$ is a general position set of $\Theta(k,\ell)$. Consequently $\gp(\Theta(k,\ell)) \geq k+1$. 
\qed

\section{Upper bounds on $\gp(G)$}
\label{sec:upper-bounds}

We say that a set of subgraphs $\{H_1,\ldots, H_k\}$ of a graph $G$ is an {\em isometric cover} of $G$ if each $H_i$, $i\in [k]$, is isometric in $G$ and $\cup_{i=1}^k V(H_i) = V(G)$. Every isometric cover of $G$ yields an upper bound on $\gp(G)$ as follows. 

\begin{theorem} [Isometric Cover Lemma]
\label{lem:iso-cover}
If $\{H_1,\ldots, H_k\}$ is an isometric cover of $G$, then 
$$\gp(G) \le \sum_{i=1}^k \gp(H_i)\,.$$ 
\end{theorem}

\proof
Let $R$ be a gp-set of $G$ and let $R_i = R\cap V(H_i)$ for $i\in [k]$. We claim that $R_i$ is a general position set of $H_i$. Suppose on the contrary that there exists vertices $x,y,z\in V(H_i)$ such that $y$ lies on some $x,z$-geodesic in $H_i$, that is, $d_{H_i}(x,z) = d_{H_i}(x,y) + d_{H_i}(y,z)$. Since $H_i$ is isometric in $G$ this implies that $d_{G}(x,z) = d_{G}(x,y) + d_{G}(y,z)$ holds, but then $R$ is not a general position set of $G$. This contradiction proves the claim. From the claim it follows that $\gp(H_i) \ge |R_i|$. We conclude that 
$$\gp(G) = |R|  = |\cup_{i=1}^k R_i| \le \sum_{i=1}^k |R_i| \le \sum_{i=1}^k \gp(H_i)\,.$$  
\qed

The {\em isometric-path number}~\cite{Fitz99, PaCh05, pan-2006} of a graph $G$, denoted by $\ip(G)$, is the minimum number of isometric paths (geodesics) required to cover the vertices of $G$. We similarly say that the {\em isometric-cycle number} of $G$, denoted by $\ic(G)$, is the minimum number of isometric cycles required to cover the vertices of $G$. Since  $\gp(P_n) \le 2$ for $n\ge 1$ and $\gp(C_n) \leq 3$ for $n\ge 3$, Isometric Cover Lemma implies:   

\begin{corollary}
\label{cor:iso-gp}
If $G$ is a graph, then 
\begin{enumerate}
\item[(i)] $\gp(G) \leq 2\ip(G)$, and
\item[(ii)] $\gp(G) \leq 3\ic(G)$.
\end{enumerate}
\end{corollary}

The bounds of Corollary~\ref{cor:iso-gp} are sharp as demonstrated by paths and complete graphs of even order for the first bound, and cycles for the second bound. 

For another upper bound we introduce the following concepts. If $v$ is a vertex of a graph $G$, then let $\ip(v,G)$ be the minimum number of isometric paths, all of them starting in $v$, that cover $V(G)$. A vertex of a graph $G$ that lies in at least one gp-set of $G$ is called a {\em gp-vertex} of $G$. Then we have: 

\begin{theorem}
\label{lem:iso-v-path-gp}
If $R$ is a general position set of a graph $G$ and $v\in R$, then
\begin{equation}
\label{math:eq2}
|R| \leq \ip(v,G) + 1
\end{equation}
In particular, if $v$ is a gp-vertex, then $\gp(G) \le \ip(v,G) + 1$. 
\end{theorem}

\proof
Let $R$ be a general position set and $v\in R$. Let $k=\ip(v,G)$. Then there exist $k$ geodesics $\{P_{vu_i}:\ u_i \in V(G), i \in [k]\}$ that covers $V(G)$. Since $R$ is a general position set, $v\in R$, and $P_{vu_i}$ is a geodesic, we have $|R \cap \left(V(P_{vu_i})\setminus \{v\} \right)|\le 1$ for $i\in [k]$. It follows that $|R| \leq k + 1 = \ip(v,G) + 1$. 

If $v$ is a gp-vertex, then consider $R$ to be a gp-set that contains $v$. By the above arguments, $\gp(G)  = |R| \le \ip(v,G) + 1$. 
\qed

If $G$ is a graph and BFS($v$) a breadth first search tree of $G$ rooted at $v$, then let $\ell(v)$ denote the number of leaves of BFS($v$). 

\begin{corollary}
\label{cor:BfsUpperBound}
If $G$ is a graph, then 
$$\gp(G) \leq 1 + \min \{\ell(v):\ v\ {\rm is\ a\ gp{\mbox -}vertex\ of}\ G\}\,.$$  
\end{corollary}

\proof
Let $v$ be a gp-vertex of $G$ and let $S$ be a gp-set containing $v$. Then $\ip(v,G) \le \ell(v)$ and hence $\gp(G) \le \ip(v,G) + 1 \le \ell(v) + 1$. Since the argument holds for any gp-vertex, the assertion follows.   
\qed

Corollary~\ref{cor:BfsUpperBound} seems in particular useful for vertex-transitive graphs because in that case it suffices to consider a single BFS tree. For a simple example consider the cycle $C_n$, $n\ge 3$. Then $\ell(v) = 2$ for any vertex $v$ of $C_n$ and hence $\gp(C_n) \le 3$ holds by Corollary~\ref{cor:BfsUpperBound}. 

To show that in Corollary~\ref{cor:BfsUpperBound} the minimum cannot be taken over all vertices, consider the following example. Let $n\ge 2$ and let $G_n$ be the graph on the vertex set $X_n \cup Y_n \cup Z_n\cup \{w\}$, where $X_n = \{x_1,\ldots, x_n\}$, $Y_n = \{y_1,\ldots, y_n\}$, and $Z_n = \{z_1,\ldots, z_n\}$. The vertices from $X_n$ induce a complete subgraph. In addition, $x_i$ is adjacent to $y_i$ and $z_i$ for $i\in [n]$, while $w$ is adjacent to all vertices from $Z_n$. Then the BFS-tree rooted in $w$ has $n$ leaves, that is, $\ell(w) = n$. On the other hand, if $u,v\in Y_n\cup Z_n$, $u\ne v$, then $d_{G_n}(u,v)\in \{2,3\}$. It follows that $Y_n\cup Z_n$ is a general position set of $G_n$, therefore $\gp(G_n) \ge 2n$.

We now turn our attention to simplicial vertices and prove: 

\begin{lemma}
\label{lem:simplicial}
If $S$ is the set of simplicial vertices of a graph $G$, then $S$ is a general position set. 
\end{lemma}

\proof
Assume on the contrary that there exist different vertices $u,v,w\in S$ such that $d_G(u,w) = d_G(u,v) + d_G(v,w)$ and let $P$ be a $u,w$-geodesic that contains $v$. Let $v'$ and $v''$ be the neighbors of $v$ on $P$, where $v'$ lies in the $u,v$-subpath of $P$ and $v''$ on the $v,w$-subpath of $P$. (Note that it is possible that $v' = u$ or $v'' = w$.) Since $v$ is a simplicial vertex, $v'v''\in E(G)$, but then $P$ is not a geodesic, a contradiction.  
\qed

\begin{figure}[ht!]
	\begin{center}
		\scalebox{0.8}{\includegraphics{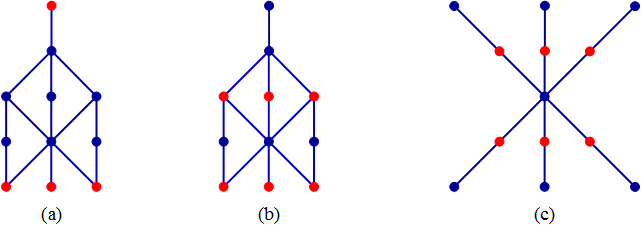}}
	\end{center}
	\caption{The red vertices form general position sets. (a) A general position set containing the unique simplicial vertex (b) A general position set without the simplicial vertex (c) A general position set with none of the six simplicial vertices.}
	\label{fig:SimplicialGPset}
\end{figure}

So simplicial vertices form general position sets. Examples in Fig.~\ref{fig:SimplicialGPset} illustrates that there is no general correlation between sets of simplicial vertices and gp-sets. However, in specific classes of graphs, the set of simplicial vertices form a gp-set. We have already noticed that this holds for complete graphs. Applying Theorem~\ref{lem:iso-v-path-gp}, in the next result we generalize this observation to all block graphs.

\begin{theorem}
\label{thm:block}
Let $S$ be the set of simplicial vertices of a block graph $G$. Then $S$ is a gp-set and hence $\gp(G) = |S|$.
\end{theorem}

\proof
Let $S$ be the set of simplicial vertices of a block graph $G$ and let $R$ be a general position set of $G$. Let $w$ be an arbitrary vertex of $R$. Since $G$ is a block graph, $w$ is is either a simplicial vertex or a cut-vertex. Hence we distinguish two cases. 

\medskip\noindent
{\bf Case 1}: $w\in S$. \\
Consider $\Psi_w = \{P_{wv}:\ v\neq w, v \in S, P_{wv}\ {\rm is\ a}\ w,v\mbox{-}{\rm geodesic}\}$. It is known~\cite{PaCh05} that $\Psi_w$ is an isometric path cover of $G$. Hence Theorem~\ref{lem:iso-v-path-gp} implies that $|R| \leq |\Psi_w| + 1 = |S|$.   
	
\medskip\noindent
{\bf Case 2}: $w\notin S$, that is, $w$ is a cut-vertex. \\
Let now $\Psi_w = \{P_{wv}:\ v \in S, P_{wv}\ {\rm is\ a}\ w,v\mbox{-}{\rm geodesic}\}$. Then again $\Psi_w$ is an isometric path cover of $G$, hence as above we have that $|R| \leq |\Psi_w| + 1 = |S| + 1$. Let now $v_1$ and $v_2$ be simplicial vertices of $G$ that are in different connected components of $G-v$. Let $P$ be the concatenation of the geodesics $P_{wv_1}$ and $P_{wv_2}$. It is easy to see that $P$ is a geodesic in $G$. Since $|R\cap V(P)| \le 2$, one of $P_{wv_1}$ and $P_{wv_2}$ intersects $R$ only in $w$. Hence, $|R| \leq (|S| + 1) - 1 = |S|$.

\medskip
We have thus proved that in both cases $|R|\le |S|$, so that $\gp(G) \le |S|$. Lemma~\ref{lem:simplicial} completes the argument. 
\qed

\begin{corollary}
\label{cor:trees}
If $L$ is the set of leaves of a tree $T$, then $\gp(T) = |L|$. 
\end{corollary}

\begin{figure}[ht!]
	\begin{center}
		\scalebox{0.25}{\includegraphics{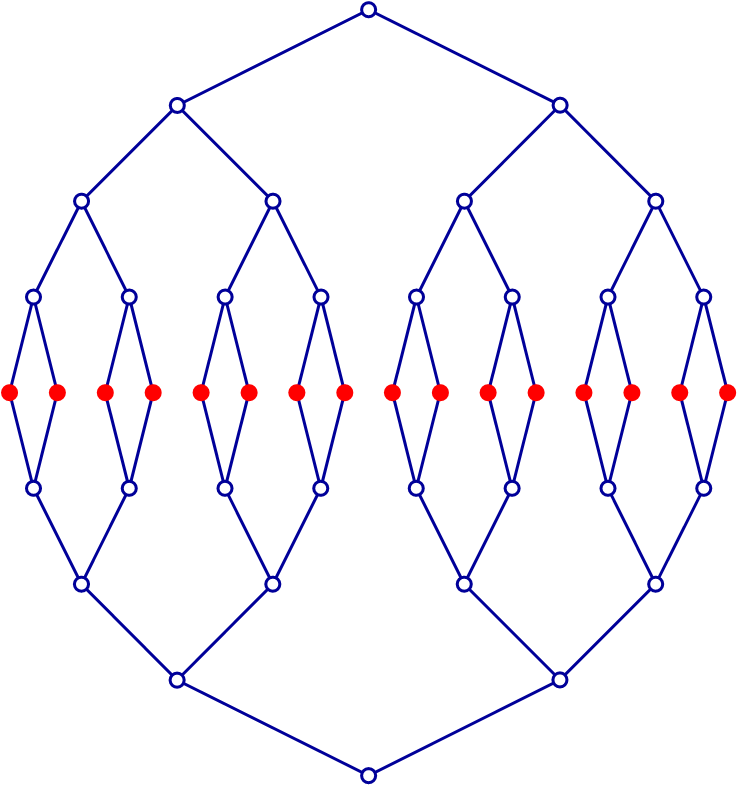}}
	\end{center}
	\caption{The glued binary tree $GT(4)$}
	\label{fig:Glued-Binary-Treea}
\end{figure}

Consider next the {\em glued binary tree} $GT(r)$, $r\ge 2$, which is obtained from two copies of the complete binary trees of depth $r$ by pairwise identifying their leaves. The construction should be clear from Fig.~\ref{fig:Glued-Binary-Treea}, where the glued binary tree $GT(4)$ is shown. The vertices obtained by identification are drawn in red, we will call then {\em quasi-leaves} of the glued binary tree.  

\begin{proposition}
\label{thm:Glued-Binary-Tree}
If $r\ge 2$, then $\gp(GT(r)) = 2^r$.
\end{proposition}

\proof
Let $R$ be a gp-set of $GT(r)$ and let $S$ be the set containing the quasi-leaves of $GT(r)$. We now consider two cases. 

\medskip\noindent
{\bf Case 1}: $R \cap S \neq \emptyset$. \\
Let $u$ be an arbitrary vertex from $R\cap S$. Then it is easy to construct geodesics $P_{uv}$ from $u$ to all other vertices $v$ of $S$ in such a way $\{P_{uv}:\ v \in S\}$ is an isometric path cover of $GT(r)$. Therefore $\gp(GT(r)) \le 1 + (|S|-1) = |S|$ by Theorem~\ref{lem:iso-v-path-gp}. 

\medskip\noindent
{\bf Case 2}: $R \cap S = \emptyset$.\\
Let now $u$ be a vertex of $R$ that is closest to a quasi-leaf among the vertices of $R$, and let $w$ be a quasi-leaf that is closest to $u$ among all quasi-leaves. Then $R' = (R \setminus \{u\}) \cup \{w\}$ is a general position set. Indeed, suppose this is not the case. Then a triple $U$ of vertices from $R'$ exists such that they lie on the same geodesic. Clearly, $w\in U$ for otherwise $R$ would not be a general position set. But then $(U\setminus \{w\})\cup \{u\}$ is a triple of vertices of $R$ lying on a common geodesic. This contradiction proves that $R'$ is a general position set. Since $|R'| = |R|$ the set $R'$ is actually a gp-set. But now we are in Case 1 and hence conclude again that $\gp(GT(r)) \le |S|$. 

\medskip
We have thus proved that $\gp(GT(r)) \leq |S|$. Since it is easy to see that $S$ is a general position set of $GT(r)$, we also have $\gp(GT(r)) \geq |S|$. We are done because $|S| = 2^r$.
\qed

\section{Lower bounds on $\gp(G)$}
\label{sec:lower-bounds}

In this section we consider lower bounds on the general position number. We already have a lower bound based on Lemma~\ref{lem:simplicial}: if $S$ is the set of simplicial vertices of a graph $G$, then $\gp(G) \geq |S|$.  

Additional lower bounds given here are in terms of the diameter of a graph and the $k$-packing number that is defined as follows. A set $S$ of vertices of a graph $G$  is a {\em $k$-packing} if $d(u,v) > k$ holds for every different $u,v\in S$. The {\it $k$-packing number} $\alpha_k(G)$ of $G$ is the cardinality of a maximum $k$-packing set~\cite{meir-1975}. For additional results on $k$-packing, see~\cite{dankelmann-2010, hartnell-1997}. Moreover, $k$-packings are the key ingredients for the concept of the $S$-packing chromatic number, see~\cite{barnaby-2017, bresar-2017, goddard-2012} and references therein. The $1$-packings are precisely independent sets and so the independence number $\alpha(G)$ is just $\alpha_1(G)$. 

A general position set need not be an independent set and vice versa. But we do have the following connection. 

\begin{proposition}
	\label{paking-gp}
	Let $G$ be a graph and $k \geq 1$. Then $\diam(G)$ $\leq$ $2k + 1$ if and only if every $k$-packing of $G$ is a general position set.
\end{proposition}

\proof
Suppose that $S$ is a $k$-packing of $G$ that is not a general position set. Then $S$ contains
vertices $x$, $y$, $z$ such that $y$ lies on an $xz$-geodesic $P_{xz}$. Since $S$ is a $k$-packing, we
have $d(x, y) \geq k + 1$ and $d(y, z) \geq k + 1$. Since $P_{xz}$ is a geodesic, it follows that $d(x, z) \geq 2k + 2$. So $\diam(G) \geq 2k + 2$.

Conversely, suppose that $\diam(G) \geq 2k + 2$. Let $x$ and $z$ be vertices with $d(x, z) = 2k+2$. In addition, let $P_{xz}$ be an $xz$-geodesic, and let $y$ be a vertex of $P_{xz}$ such that $d(x, y) = d(y, z) = k + 1$. Then $\{x, y, z\}$ is a $k$-packing that is not a general position set.
\qed

Proposition~\ref{paking-gp} provides a lower bound on $\gp$-sets of a graph.
\begin{corollary}
\label{cor:packing-bound}
If $G$ is a graph with $\diam(G)$ $\leq$ $2k + 1$, then $\gp(G) \geq \alpha_k(G)$.
\end{corollary}

Since 1-packing sets of a graph are precisely its independent sets, Proposition~\ref{paking-gp} for $k = 1$ asserts the following corollary: 
\begin{corollary}
	\label{prop:stable-gp}
	Let $G$ be a graph and $k \geq 1$. Then $\diam(G)$ $\leq$ $3$ if and only if every independent set of $G$ is a general position set.
\end{corollary}
In general, however, there is no connection between the independence number $\alpha(G)$ of $G$ and $\gp(G)$. For instance, $\gp(K_n)=n$ and $\alpha(K_n) = 1$, while on the other hand, $\gp(P_n)$ = 2 and $\alpha(P_n) = \lceil n/2 \rceil$. 

Another lower bound on $\gp(G)$ involves the distance between the edges of a graph which is defined as follows. If $e=uv$ and $f=xy$ are edges of a graph $G$, then $d(e,f)$ = $\min\{d(u,x), d(u,y), d(v,x), d(v,y)\}$. 

\begin{proposition}
	\label{prp:set-of-edges}
	Let $G$ be a graph with $\diam(G) = k \ge 2$. If $F$ is a set of edges of $G$ such that $d(e,f) = k$ for every $e,f\in F$, $e\ne f$, then $\gp(G)\ge 2|F|$.  
\end{proposition}
\proof
We claim that the set $S$ consisting of the end-vertices of the edges from $F$ is a general position set. If $x\in S$, then let $f_x$ be the edge of $F$ containing $x$. Let $x,y,z$ be an arbitrary triple of vertices from $S$ and suppose that $y$ lies on a $x,z$-geodesic $P$. Clearly, $f_x\ne f_z$. Since $d(f_x,f_z) = k$ and $\diam(G) = k$, we must necessarily have $d(x,z) = k$. Suppose without loss of generality that $f_y\ne f_x$. But then, as $P$ is a geodesic, we have $d(y,x)\le k-1$ and hence $d(f_x,f_y)\le k-1$, a contradiction. 
\qed

\begin{figure}[ht!]
	\begin{center}
		\scalebox{0.25}{\includegraphics{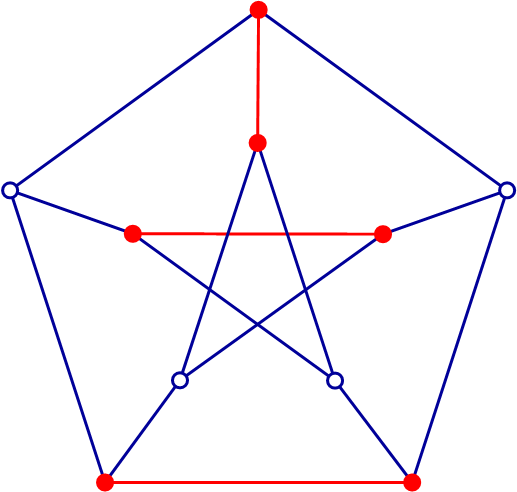}}
	\end{center}
	\caption{$\gp(P)\ge 6$}
	\label{fig:petersen}
\end{figure}

As an application of the above proposition, consider the Petersen graph $P$. On Fig.~\ref{fig:petersen} three edges of $P$ pairwise at distance $2$ are shown, hence $\gp(P) \ge 6$ by Proposition~\ref{prp:set-of-edges}. Since $V(P)$ can be covered with two disjoint isometric cycles,  $\gp(P) \ge 6$ by  Corollary~\ref{cor:iso-gp}(ii). Thus $\gp(P) = 6$.

Additional examples demonstrating sharpness of Proposition~\ref{prp:set-of-edges} where $k$ is large, can be constructed as follows. Start with the star $K_{1,n}$ and subdivide each edge of it the same number of times. Then to each of the $n$ leaves attach a private triangle by identifying a vertex of the triangle with the leaf. The $n$ edges of these triangles whose end-vertices are of degree $2$ are edges that satisfy the assumption(s) of the proposition.  

\section{Computational complexity of the problem}
\label{sec:NPComp}
The (graph) general position problem is the following:

\begin{center}
	\fbox{\parbox{0.85\linewidth}{\noindent
			{\sc General Position Problem}\\[.8ex]
			\begin{tabular*}{.93\textwidth}{rl}
				{\em Input:} & A graph $G$, and an integer $k$.\\
				{\em Question:} & Is $\gp(G)\ge k$?
			\end{tabular*}
		}}
	\end{center}
	
	The general position subset selection problem from discrete geometry which is a main  motivation of this paper has been proved as NP-hard~\cite{FrKa16, PaWo13}. We next prove a parallel result for the {\sc General Position Problem}. 
	
	\begin{theorem}
		\label{thm:NP-complete}
		{\sc General Position Problem} is NP-complete.
	\end{theorem}
	
	\proof
	Note first that the {\sc General Position Problem} is in NP. A set $S$ of vertices of a graph $G$ is a general position set of $G$ if and only if for each pair of vertices $x$ and $z$ of $S$, $d(x,z) \neq d(x,y)+d(y,z)$ for every $y$ in $S$. This task can be done in polynomial time.  In the rest of the proof, we give a reduction of the NP-complete {\sc Maximum Independent Set Problem}, to the {\sc General Position Problem}. The former problem is one of the classical NP-complete problems~\cite{GaJo79}. 
	
	Given a graph $G = (V,E)$, we construct a graph $\widetilde{G} = (\widetilde{V},\widetilde{E})$  as follows. Its vertex set is $\widetilde{V} = V\cup V'\cup V''$, where $V'=\{v':\ v \in V\}$  and $V'' = \{v'':\ v \in V\}$. The set of edges is $\widetilde{E}=E\cup E'\cup E''\cup E'''$, where $E'$ is the set of all possible edges between the vertices of $V'$, while  $E''=\{v,v':\ v \in V\}$ and $E'''=\{v'v'':\ v\in V\}$. The graph $\widetilde{G}$ can be thus considered as composed of three parts: the original graph $G$, the complete graph induced by $V'$, and the independence set induced by $V''$. These three parts are connected by the matching $E''$ between $V$ and $V'$ and the matching $E'''$ between $V'$ and $V''$. 
	
	We first claim that $X\subseteq V$ is an independent set of $G$ if and only of $X\cup V_2$ is a general position set of $\widetilde{G}$. Suppose first that $X\subseteq V$ is an independent set of $G$. Then, clearly, $X\cup V''$ is an independent set of $\widetilde{G}$. Since $\diam(\widetilde{G}) = 3$, Corollary~\ref{prop:stable-gp} implies that $X\cup V''$ is a general position set of $\widetilde{G}$. Conversely, assume that $X$ is not independent and let $x,y\in X$ be adjacent vertices. Then the path $xyy'y''$ is a geodesic which in turn implies that $X\cup V''$ is not a general position set of $\widetilde{G}$. 
	
	We next claim that $\alpha(G)\ge k$ if and only if $\gp(\widetilde{G}) \ge k + |V|$.  It suffices to show that if $S$ is a general position set of $\widetilde{G}$, then there exists a general position set $\widetilde{S}$ of $\widetilde{G}$ such that $\widetilde{S} = X\cup V''$, where $X$ is an independent set of $G$ and $|\widetilde{S}| \ge |S|$. For any two vertices $x$ and $y$ of $V$ and its corresponding vertices $x'$ and $y'$ of $V'$ and $x''$ and $y''$ of $V''$,  $x''x'y'y''$ is a geodesic in $\widetilde{G}$. For some $u \in V$, if both $u'$ and $u''$ are in $S$, then no other vertices $v'' \neq u''$ will be in $S$. This will contradict the maximality of $\gp$-set of  $\widetilde{G}$ when $|V| \geq 3$.  If $u' \in S$ and $u'' \notin S$, then consider $\widetilde{S} = S \cup \{u''\} \setminus \{u'\}$. It concludes that given a general position set $S$ of $\widetilde{G}$, there exists a general position $\widetilde{S}$ of $\widetilde{G}$ such that $\widetilde{S} \cap V' = \emptyset$. From here the claim follows which also completes the argument.  
	\qed
	
\section{Concluding remarks}
\label{sec:concluding remarks}

In this paper, a new graph combinatorial problem  is introduced. Even though the discrete  geometry version of the general position problem has been well-studied over the years, its graph theory version has not been investigated to the best of our knowledge. In this paper we prove that the graph theory general position problem is NP-complete. Connections between the general position sets, packings, and simplicial vertices are studied. Using isometric path covers and isometric cycle covers sharp lower and upper bounds on the general position number were derived. With the aid of these these bounds we have solved the general position problem for block graphs, theta graphs, and glued-binary trees.

Since this problem is new, the research topic is wide open. It will be interesting to study the general position problem for important interconnection networks such as butterfly, grid-like architectures, etc., as well as for Cayley graphs in general and their subclasses such as torus graphs and hypercubes in particular. In the same way, the complexity status of the general position problem for important classes of graphs such as bipartite graphs, chordal graphs, and planar graphs is open. 

\section*{Acknowledgments}
	
This work was supported and funded by Kuwait University, Research Project No.\ (QI 01/16).


\end{document}